\documentclass[leqno]{article}
\usepackage{verbatim, amsmath, amscd, amssymb}
\usepackage[all]{xy}

\usepackage[frenchb]{babel}

\newcommand{\rC}{\mathrm{C}}

\newcommand{\rG}{\mathrm{G}}

\newcommand{\rR}{\mathrm{R}}
\newcommand{\rs}{\mathrm{s}}

\newcommand{\bbB}{\mathbb B}

\newcommand{\bbF}{\mathbb F}

\newcommand{\bbN}{\mathbb N}
\newcommand{\bbQ}{\mathbb Q}

\newcommand{\bbZ}{\mathbb Z}

\newcommand{\ch}{\mathrm{ch}\,}

\newcommand{\Dist}{\mathrm{Dist}}

\newcommand{\Ext}{\mathrm{Ext}}

\newcommand{\id}{\mathrm{id}}

\newcommand{\ind}{\mathrm{ind}}
\newcommand{\coind}{\mathrm{coind}}

\newcommand{\SL}{\mathrm{SL}}

\newcommand{\St}{\mathrm{St}}

\newcommand{\Mod}{\mathbf{Mod}}
\newcommand{\Sch}{\mathbf{Sch}}

\newcommand{\lbr}{\begin{bmatrix}}
\newcommand{\rbr}{\end{bmatrix}}
\newcommand{\for}{\bigcirc\kern-2.6ex \because}
\newcommand{\forb}{\bigcirc\kern-2.8ex \because}
\newcommand{\forbb}{\bigcirc\kern-3.0ex \because}
\newcommand{\forbbb}{\bigcirc\kern-3.1ex \because}

\newcommand\pf{\noindent {\it Preuve.   }}
\newtheorem{thm}{Th\'eor\`eme.}

\newtheorem{prop}{Proposition.}

\newtheorem{lem}{Lemme.}

\newtheorem{cor}{Corollaire.}

\begin{document}
\large
\title{
{\bf 
Contraction par Frobenius
et modules de Steinberg 
}
\thanks{Le second auteur a b\'en\'efici\'e lors de ce travail  d'un soutien  JSPS Grants in Aid for Scientific Research 15K04789
}
\author{
Michel Gros
\\
CNRS UMR 6625, IRMAR, 
\\
Universit\'e de Rennes 1,
Campus de Beaulieu, 
\\
35042 Rennes cedex, France
\\
michel.gros@univ-rennes1.fr
\and
K\textsc{aneda} Masaharu
\\
Osaka City University
\\
Department of Mathematics
\\
Osaka 558-8585, Japan
\\
kaneda@sci.osaka-cu.ac.jp
}
}
\date{}
\maketitle

\begin{abstract}

Soit  $G$ un groupe r\'eductif d\'efini sur un corps alg\'ebriquement clos de caract\'eristique positive. Nous montrons que le foncteur \emph{contraction par Frobenius} de la cat\'egorie des $G$-modules est adjoint \`a droite de celui de tensorisation deux fois par le module de Steinberg du tordu de Frobenius du $G$-module de d\'epart. Il s'ensuit en particulier que le foncteur de contraction par Frobenius pr\'eserve le caract\`ere injectif et l'existence de bonne filtration mais pas la semi-simplicit\'e.\\

For a reductive group $G$ defined over an algebraically closed field of positive characteristic,
we show that the    \emph{Frobenius contraction} functor of $G$-modules 
is right adjoint to the Frobenius twist of the modules tensored with the Steinberg module twice.
It follows that the Frobenius contraction functor preserves
injectivity, good filtrations, but not semisiplicity.

\end{abstract}

Soit $G$ un groupe alg\'ebrique r\'eductif d\'efini sur un corps alg\'ebriquement clos $\Bbbk$ de caract\'eristique $p>0$.  
Le morphisme induit sur l'alg\`ebre des distributions $\Dist(G)$ de $G$ par l'endomorphisme de Frobenius de $G$ 
admet un scindage \cite{GK11}, \cite{GK}.
Ce dernier permet de d\'efinir une nouvelle op\'eration sur les $G$-modules, qu'\`a la suite de Littelmann qui l'avait 
d\'efini auparavant seulement pour les modules de Weyl duaux, 
nous appellons {\emph{contraction par Frobenius}}.
Si $M$ est un $G$-module,
la contraction par  Frobenius  permet de munir la somme directe des sous-espaces de $M$ de poids divisible
par $p$ d'une structure de $G$-module. Cette construction  g\'en\'eralise le processus de ``d\'etorsion'' par Frobenius sur $M$ 
not\'e $M^{[-1]}$
dans  \cite[II.3.16]{J}, lorsque l'on part d'une situation o\`u le noyau de Frobenius de  $G$ agit trivialement sur $M$.

Nous pr\'esentons tout d'abord dans cet article la caract\'erisation,  d\'egag\'ee par 
Stephen Donkin, du foncteur  contraction par Frobenius comme adjoint 
\`a droite de celui de double tensorisation par le module de Steinberg du tordu par Frobenius (2.4, Cor.).  Ceci nous permet ensuite d'\'etablir   inconditionnellement quelques propri\'et\'es du foncteur  contraction par Frobenius  : pr\'eservation du caract\`ere injectif et de l'existence d'une bonne filtration pour les $G$-modules (3.1, Th.) mais pas de la semi-simplicit\'e en g\'en\'eral (1.4). Enfin, nous d\'emontrons des r\'esultats similaires ($4^\circ$) pour les $G_r$-modules,  $G_r$ d\'esignant  le $r$-i\`eme noyau de Frobenius de   $G$.
 
Cet  article tire son origine   d'une question pos\'ee \`a S. Donkin par  le second auteur (M.K.)  lors de sa visite \`a l'Institut  
Mittag-Leffler en Mai 2015. Il remercie cet Institut de lui avoir procur\'e cette opportunit\'e 
ainsi que   Donna Testerman de l'avoir invit\'e \`a exposer une partie de ce travail \`a l'Institut Bernouilli en Ao\^ut 2016. Les auteurs remercient chaleureusement S. Donkin de leur avoir permis d'inclure certains de ses r\'esultats  et assument la pleine  responsabilit\'e  d'\'eventuelles erreurs.\\

\setcounter{equation}{0}
\begin{center}
$1^\circ$
{\bf 
Scindage du Frobenius de $\Dist(G)$}
\end{center}

Soient  $G$ comme ci-dessus et  $F$ son  morphisme  de Frobenius.
Par commodit\'e, nous consid\`ererons une   $\bbF_p$-forme $G_{\bbF_p}$
de $G$ et de son endomorphisme de Frobenius g\'eom\'etrique d\'efini par  
$?^p\otimes_{\bbF_p}\Bbbk$ sur l'alg\`ebre des coordonn\'ees
$\Bbbk[G]=\bbF_p[G_{\bbF_p}]\otimes_{\bbF_p}\Bbbk$ de $G$
\cite[I.9.2]{J}.
On supposera \'egalement $G$ simplement connexe et semi-simple.
Soient $B$ un sous-groupe de Borel de $G$,
$T$ un tore maximal de $B$ tous deux scind\'es sur  $\bbF_p$,
$R$ le syst\`eme de racines de $G$ relativement \`a $T$, $R^+$ la partie positive de $R$ pour 
laquelle les racines de  $B$ sont $-R^+$, et $R^\rs$ l'ensemble des racines simples de $R^+$.
Pour    tout  sous-sch\'ema en groupes $H$  de  $G$, 
 $\Dist(H)$ d\'esignera l'alg\`ebre des distributions sur 
 $H$, et
$H_1$ le noyau de Frobenius de  $H$.

\noindent
(1.1)
Comme $T$ est ab\'elien, l'unit\'e de $\Dist(T_1)$ admet une unique d\'ecomposition en idempotents primitifs orthogonaux    \cite[\S1. Thm. 4.6]{NT}.
Parmi ces derniers, il en existe  un unique  
$\mu_0$ \`a ne pas \^etre  annul\'e par la co-unit\'e de $\Dist(T_1)$ et c'est aussi l'unique idempotent 
\`a ne pas \^etre annul\'e par $\Dist(F)$ ;
$\Dist(F)(\mu_0)=1$.
Soient $x_\alpha$, $\alpha\in R$, 
$h_\beta$, $\beta\in R^\rs$, une base de Chevalley de l'alg\`ebre de Lie de   $G$ correspondant \`a  $T$ et aux sous-groupes radiciels de  $G$.
Pour
$n\in\bbN$, nous noterons
$x_\alpha^{(n)}$ la $n$-i\`eme  puissance divis\'ee  de $x_\alpha$.
Il existe un  homomorphisme de $\Bbbk$-alg\`ebres
$\phi:\Dist(G)\to\mu_0\Dist(G)\mu_0$ tel que
$x_{\pm\alpha}^{(n)}\mapsto
\mu_0x_{\pm\alpha}^{(pn)}\mu_0$, 
$\binom{h_\alpha}{m}\mapsto
\mu_0\binom{h_\alpha}{pm}\mu_0$ 
pour tous  $\alpha\in R^\rs$, $n\in\bbN$,
et v\'erifiant
$\Dist(F)\circ\phi=\id_{\Dist(G)}$
\cite{GK11}, \cite{GK}.
Nous appellerons $\phi$ le {\emph{scindage du Frobenius}} sur
$\Dist(G)$.

Etant donn\'e un $G$-module $M$, en faisant agir $\Dist(G)$ sur $M$ via $\phi$, i.e. sur  $\mu_0M$, on obtient ainsi un nouveau   $\Dist(G)$-module.
Comme  $M$ est localement fini, il en est de m\^eme  de  $\mu_0M$, et l'on obtient donc (\cite[II.1.20]{J}) une structure de 
 $G$-module sur $\mu_0M$
\cite[II.1.20]{J}, qu'on appellera la {\emph{contraction par Frobenius}} de $M$ et qu'on notera  $M^\phi$.
L'action de  $\Dist(G)$ sur  $M^\phi$ sera not\'ee par  $\bullet$:
$\mu\bullet m=\phi(\mu)m$,
$\mu\in\Dist(G)$, $m\in M^\phi$.

\setcounter{equation}{0}
\noindent
(1.2)
Rappelons maintenant quelques propri\'et\'es basiques de cette contraction par Frobenius.
Comme $\mu_0$ est un  idempotent, 

(i)
le foncteur de contraction par Frobenius est   exact.

\noindent
Soit
$\Lambda$ le groupe des caract\`eres de  $T$. Pour tout $\lambda\in\Lambda$, nous noterons encore   abusivement 
$\lambda$ l'application  $\Dist(\lambda):\Dist(T)\to\Bbbk$.
Soit
$M$ un $G$-module.
Si $M_\lambda=\{m\in M|
tm=\lambda(t)m,\  \forall t\in T\}
=
\{m\in M|
\mu m=\lambda(\mu)m,\ \forall \mu\in \Dist(T)\}$ d\'esigne le sous-espace de poids  $\lambda$ de
$M$, on a :

(ii) $\mu_0M=\coprod_{\lambda\in\Lambda}M_{p\lambda}$
avec $M_{p\lambda}$
l'espace de poids 
$\lambda$ de $M^\phi$;
consid\'erant
$m\in M_{p\lambda}$
comme un \'el\'ement de
$M^\phi$,
on a 
$\binom{h_\beta}{n}\bullet m=
\binom{\langle\lambda,\beta^\vee\rangle}{n}m$
pour tous  $\beta\in R^\rs$
et $n\in\bbN$.
En particulier,
si
$\sum_{\lambda\in\Lambda}M_{p\lambda}=0$,
la contraction par Frobenius annule $M$.

(iii)
Si l'on fait agir  $G$ sur $M$ via le morphisme de Frobenius de $G$, on obtient un nouveau $G$-module appel\'e le tordu par Frobenius de  $M$ et not\'e $M^{[1]}$. 
La contraction par Frobenius de ce dernier redonne alors
$M$:
$
(M^{[1]})^\phi\simeq M$.
En particulier, le foncteur de torsion par Frobenius  
$M\mapsto M^{[1]}$ est un endofoncteur pleinement fid\`ele de la cat\'egorie des 
$G$-modules.

(iv)
Si
$G_1$ agit trivialement sur $M$,
la $G$-action sur $M$
se factorise par le quotient
$G\to G/G_1$,
que l'on peut identifier au morphisme de Frobenius  gr\^ace au diagramme commutatif :

\[
\xymatrix{
G\ar[rr]^-{F}
\ar[d]
&&
G.
\\
G/G_1
\ar@{.>}[rru]_\sim
}
\]
Il s'ensuit qu'il existe une structure de  
$G$-module, not\'ee $M^{[-1]}$
dans \cite[II.3.16]{J}, sur le  $\Bbbk$-espace vectoriel
$M
$
d\'etordant la torsion par Frobenius,  et que l'on retrouve le $G$-module $M$ de d\'epart 
en tordant par Frobenius  :
$(M^{[-1]})^{[1]}\simeq
M$. Le $G$-module
$M^{[-1]}$
n'est donc alors rien d'autre dans ce cas que la contraction par Frobenius 
$M^{\phi}$ de $M$.

(v)
La contraction par Frobenius commute avec la formation du dual :
$(M^\phi)^*\simeq(M^*)^\phi$
pour tout $M$ de dimension finie
\cite[5.5]{GK}, \cite[Prop. 6.6]{GK'}.

(vi)
Si $V$ est un   $G$-module,
$\Dist(G)$ agit sur $M\otimes V$ via la comultiplication $\Delta$ de
$\Dist(G)$;
si $\Delta(x)=\sum_ix_i\otimes
y_i$ et si $x\in\Dist(G)$, alors
$x\bullet(m\otimes v)=
\sum_ix_im\otimes
y_iv$,
$m\in M$, $v\in V$.
Bien que la contraction par  Frobenius   ne commute pas avec la comultiplication, on a :
\cite[Lem. 6.8]{GK'}
\[
(M\otimes V^{[1]})^\phi\simeq
M^\phi\otimes V,
\quad
(M^{[1]}\otimes
V)^\phi\simeq
M\otimes V^\phi.
\]

\setcounter{equation}{0}
\noindent
(1.3)
Soit
$\Lambda^+=\{\lambda\in\Lambda|\langle\lambda,\alpha^\vee\rangle\geq0,\ \forall \alpha\in R^+\}$ l'ensemble des poids dominants.
Les $G$-modules simples sont   param\'etris\'es par leur plus haut poids dominant ; notons 
alors $L(\lambda)$ le $G$-module simple de plus haut poids   $\lambda$.
Si
$\Lambda_1=\{\mu\in\Lambda^+|\langle\mu,\alpha^\vee\rangle<p,\ \forall\alpha\in R^\rs\}$,
le th\'eor\`eme de Steinberg de d\'ecomposition en  produit tensoriel 
\cite[II.3.17]{J}
donne
un isomorphisme de $G$-modules
$L(\lambda)\simeq
L(\lambda^0)\otimes
L(\lambda^1)^{[1]}\otimes
\dots\otimes L(\lambda^r)^{[r]}$,
$\lambda^i\in\Lambda_1$
avec
$\lambda=\sum_{i=0}^rp^i\lambda^i$,
et $?^{[j]}$ la torsion par $F^j$.
On a donc
\[
L(\lambda)^\phi\simeq
L(\lambda^0)^\phi\otimes
L(\lambda^1)\otimes
L(\lambda^2)^{[1]}\dots\otimes L(\lambda^r)^{[r-1]}.
\]
Une estimation grossi\`ere des poids donne alors 
\begin{prop} 
Soit $h$ le nombre de  Coxeter de
$G$. Supposons que $p\geq2(h-1)$.
Un $G$-module simple de plus haut poids appartenant à $\Lambda_1$
se contracte par Frobenius en un $G$-module semi-simple
si ce  dernier  est non nul.
\end{prop}

\pf Soit
$\lambda\in \Lambda_1$, et
soit
$p\mu$,
$\mu\in\Lambda^+$, un poids de $L(\lambda)$.
On a alors
$p\mu=\lambda-\gamma$
pour un certain
$\gamma\in\sum_{\alpha\in R^\rs}\bbN\alpha$.
Si $\alpha_0^\vee$ est la plus grande coracine de $R^\vee$,
\begin{align*}
\langle\mu+\rho,\alpha_0^\vee\rangle
&=
\langle\frac{1}{p}(\lambda-\gamma)+\rho,\alpha_0^\vee\rangle
=
\frac{1}{p}\langle\lambda-\gamma,\alpha_0^\vee\rangle+(h-1)
\\
&\leq
\frac{1}{p}\langle\lambda,\alpha_0^\vee\rangle+(h-1)
\leq
\frac{1}{p}\langle(p-1)\rho,\alpha_0^\vee\rangle+(h-1)
\\
&=
\frac{1}{p}(p-1)(h-1)+(h-1)
=
(h-1)(2-\frac{1}{p})
<2(h-1)\leq p.
\end{align*}
L'assertion est alors une cons\'equence du 
{\it linkage principle}
\cite[II.6.17]{J}.

\setcounter{equation}{0}
\noindent
(1.4)
La conclusion de la proposition est inexacte sans hypoth\`ese de caract\'eristique ou de petitesse du plus haut poids comme on le montrera en   (3.4).
D'autre part,
prendre les points fixes sous  $G_1$ d'un module semi-simple pr\'eserve, si cela ne l'annule pas, la semi-simplicit\'e.

\setcounter{equation}{0}
\begin{center}
$2^\circ$
{\bf 
Caract\'erisation par les modules de Steinberg }

\end{center}

Pour un $G$-module $M$
et pour un sous-sch\'ema en groupes $H$ de $G$
on note $M^H$ l'ensemble des \'el\'ements $H$-invariants de $M$.
Comme $G_1$ est normal dans  $G$,
$M^{G_1}$ est  $G$-sous-module de $M$.
Soit $\St=L((p-1)\rho)$,
$\rho=\frac{1}{2}\sum_{\alpha\in R^+}\alpha$,
le module de Steinberg de $G$. 
Il est auto-dual.
\\

\setcounter{equation}{0}
\noindent
(2.1)
{\bf Th\'eor\`eme
(Donkin).}
{\it
Pour tout $G$-module $M$, il existe  un isomorphisme de
$G$-modules
\[
M^\phi\simeq
\{(\St\otimes\St\otimes M)^{G_1}\}^{[-1]}.
\] 
}

\setcounter{equation}{0}
\noindent
(2.2)
Nous aurons tout d'abord besoin du lemme suivant.
Rappelons que, comme  $G_1$-module, on a $\St=\Dist(G_1)\otimes_{\Dist(B_1^+)}(p-1)\rho$
avec $B^+$ le sous-groupe de Borel oppos\'e de  $B$
\cite[II.3.18]{J}.
Fixons  un vecteur de plus haut (resp. bas) poids  
$v_+$
(resp. $v_-$) de
$\St$.
Rappelons aussi que la $\Dist(T_1)$-repr\'esentation
$\mu_0=\Dist(T_1)\mu_0$ est   triviale.

\noindent
{\bf Lemme (Donkin).}
{\it
Il existe un  $G_1$-isomorphisme
$\theta:
\Dist(G_1)\otimes_{\Dist(T_1)}\mu_0\to\St\otimes\St$
tel que  $x\otimes1
\mapsto x(v_+\otimes v_-)$.
}

\pf
On a des isomorphismes  $\Bbbk$-lin\'eaires
\begin{align*}
\Dist(G_1)\Mod(\Dist(G_1)\otimes_{\Dist(T_1)}\mu_0,
\St\otimes\St)
&\simeq
\Dist(T_1)\Mod(\mu_0,
\St\otimes\St)
\\
&\simeq
\mu_0(\St\otimes\St)=(\St\otimes\St)^{T_1}.
\end{align*}
Le choix de $v_+\otimes v_-\in(\St\otimes\St)^{T_1}$
d\'efinit donc une application    $G_1$-lin\'eaire
$\theta$ comme dans l'\'enonc\'e du lemme.
Montrons qu'elle  est injective, et donc bijective pour une raison de dimension.

Soit $U$ (resp. $U^+$) le radical unipotent de
$B$
(resp. $B^+$),
il existe des isomorphismes $\Bbbk$-lin\'eaires
$\Dist(G_1)\otimes_{\Dist(T_1)}\mu_0\simeq
\Dist(U_1)\otimes\Dist(U_1^+)$
via
$xy\otimes1\gets\!\shortmid x\otimes y$,
et 
$\Dist(U_1)\simeq\St\simeq
\Dist(U_1^+)$ via $x\mapsto xv_+$
et $yv_-\gets\!\shortmid
y$.
\'Ecrivons
$R^+=\{\alpha_1,\dots,\alpha_N\}$
avec $N=|R^+|$
et munissons la base de  PBW de $\Dist(U_1)$ de l'ordre lexicographique tel que   
$\prod_{i=1}^Nx_{-\alpha_i}^{(n_i)}
\prec
 \prod_{i=1}^Nx_{-\alpha_i}^{(m_i)}$
si et seulement si, soit   $n_N<m_N$, soit
il existe  $j\in[1,N]$ tel que
 $n_i=m_i$ pour tout $i\in[j,N]$ et $n_{j-1}<m_{j-1}$.
Chaque
$x_{\alpha}^{(n)}$, $\alpha\in R$, $n\in\bbN$,
agit sur
$\St\otimes\St$ par
$\Delta(x_{\alpha}^{(n)})=\sum_{i=0}^nx_{\alpha}^{(i)}\otimes
x_{\alpha}^{(n-i)}$
\cite[I.7.8]{J}.
Alors
\[
\theta(x\otimes y\otimes1)=
xy(v_+\otimes v_-)
=
x(v_+\otimes yv_-)
\]
pour tous 
$x\in\Dist(U_1)$ et $y\in\Dist(U_1^+)$
car
$v_+$ est de poids maximum. 
Supposons qu'il existe
 $\sum_{(n_i)} c_{(n_i)}(\prod_{i=1}^N
 x_{-\alpha_i}^{(n_i)})\otimes y$  
 avec  $c_{(n_i)}
 \in\Bbbk\setminus0$,
 $y\in\Dist(U_1^+)\setminus0$,
tels que
 $\{\sum c_{(n_i)}(\prod_i
 x_{-\alpha_i}^{(n_i)})\otimes y\}(v_+\otimes v_-)=0$.
Si $(m_i)$ fournit l'exposant conduisant au plus grand \'el\'ement (pour l'ordre lexicographique) parmi les
$\prod
 x_{-\alpha_i}^{(n_i)}$,
 alors
\[
\{\sum c_{(n_i)}(\prod_i
x_{-\alpha_i}^{(n_i)})\otimes y\}(v_+\otimes v_-)\in
c_{(m_i)}(\prod_i
x_{-\alpha_i}^{(m_i)}v_+)\otimes yv_-+
\sum 
\Bbbk
x'v_+\otimes x''yv_+,
\]
expression dans laquelle la seconde somme \`a droite porte sur certains \'el\'ements 
$x',x''$ de la base de PBW de $\Dist(U_1)$ avec
$x'\prec\prod
x_{-\alpha_i}^{(m_i)}$, ce qui contredit l'hypoth\`ese que $c_{(m_i)}\ne0$.
\quad

\setcounter{equation}{0}
\noindent
(2.3)
Venons-en \`a la preuve du th\'eor\`eme.
Le lemme montre que l'on peut munir  
$\mu_0M
$
d'une structure de $G$-module via des   isomorphismes $\Bbbk$-lin\'eaires
\begin{align*}
\mu_0M
&\simeq
\Dist(T_1)(\mu_0, M)
\simeq
\Dist(G_1)\Mod(\Dist(G_1)\otimes_{\Dist(T_1)}\mu_0, M)
\\
&\simeq
G_1\Mod(\St\otimes\St, M),
\end{align*}
se d\'ecrivant comme
\begin{equation}
\mu_0M\ni
m\longmapsto f_m\in
G_1\Mod(\St\otimes\St, M) \quad\text{tel que}\quad
f_m(v_+\otimes v_-)=m.
\end{equation}
Pour $x\in\Dist(G_1)$,
$f_m(x(v_+\otimes v_-))=x m$
avec
$x$
agissant sur  $m$ dans le membre de droite
suivant la structure de $G$-module sur $M$.

Comme
$(\St\otimes\St\otimes M)^{G_1}\simeq
G_1\Mod(\St\otimes\St, M)$
par auto-dualit\'e de 
$\St$,
on va montrer que ce dernier est isomorphe \`a  
$(M^\phi)^{[1]}$
comme $G$-module.
On est ainsi r\'eduit \`a montrer  
que l'application 
\[
\eta:(M^\phi)^{[1]}\to
G_1\Mod(\St\otimes\St, M)
\quad\text{via}\quad
m\mapsto f_m
\]
est
$G$-lin\'eaire.
Comme
$\Dist(G)$ est engendr\'e par les
$x_{\pm\alpha}^{(n)},
\alpha\in R^\rs, n\in\bbN$
\cite[Satz I.7]{J73},
il suffit de v\'erifier que 
\begin{equation}
\eta(x_{\pm\alpha}^{(n)}\bullet m)=
x_{\pm\alpha}^{(n)}\eta(m),
\end{equation}
expression dans laquelle
$\bullet$ dans le membre de gauche d\'esigne l'action de 
$x_{\alpha}^{(n)}$ sur
$(M^\phi)^{[1]}$.
Rappelons que l'action de  
$\Dist(G)$ sur $G_1\Mod(M,M')$
pour des  $G$-modules $M$ et $M'$ est telle que pour $f\in
G_1\Mod(M,M')$
et $\beta\in R$, on ait 
$x_{\beta}^{(n)}f=\sum_{i=0}^nx_\beta^{(n-i)}f((-1)^ix_\beta^{(i)}?)$
\cite[I.7.8, 11]{J}
gr\^ace \`a la bijection
$(M^*\otimes M')^{G_1}\to
G_1\Mod(M,M')$
via
$\gamma\otimes m'\mapsto
\gamma(?)m'$.

On a donc 
\begin{align*}
(x_{\alpha}^{(n)} f_m)
&
(v_+\otimes v_-)
=
\sum_{b\in[0,n]}x_\alpha^{(n-b)} f_m((-1)^bx_\alpha^{(b)}(v_+\otimes v_-))
\\
&=
\sum_{b\in[0,n]}x_\alpha^{(n-b)}f_m((-1)^b\sum_{s\in[0,b]}(x_\alpha^{(s)}\otimes x_\alpha^{(b-s)})(v_+\otimes v_-))
\\
&=
\sum_{b\in[0,n]}(-1)^bx_\alpha^{(n-b)} f_m(v_+\otimes x_\alpha^{(b)}v_-)
\\
&\hspace{3cm}
\quad\text{car $v_+$ est un vecteur de plus haut poids de   $\St$}
\\
&=
\sum_{b\in[0,p[}(-1)^bx_\alpha^{(n-b)}
f_m(v_+\otimes x_\alpha^{(b)}v_-)
\\
&\hspace{3cm}
\quad\text{car $-(p-1)\rho+b\alpha$ n'est pas un poids de $\St$ pour $b\geq p$}
\\
&=
\sum_{b\in[0,p[}(-1)^bx_\alpha^{(n-b)} f_m(x_\alpha^{(b)}(v_+\otimes v_-))
\\
&\hspace{3cm}
\quad\text{car,  \`a nouveau, $v_+$ est  un vecteur de plus haut poids de 
$\St$ }
\\
&=
\sum_{b\in[0,p[}(-1)^bx_\alpha^{(n-b)}x_\alpha^{(b)} f_m(v_+\otimes v_-)
\quad\text{car $ f_m$ est $G_1$-lin\'eaire}
\\
&=
\sum_{b\in[0,p[}
(-1)^b\binom{n}{b}
x_\alpha^{(n)}m
\\
&=
\sum_{b\in[0,p[}
(-1)^b\binom{n_0}{b}
x_\alpha^{(n)}m
\quad\text{en \'ecrivant 
$n=n_0+pn_1$ avec
$n_0\in[0,p[$ et $n_1\in\bbN$}
\\
&=
\sum_{b\in[0,n_0]}
(-1)^b\binom{n_0}{b}
x_\alpha^{(n)}m
=
(1-1)^{n_0}
x_\alpha^{(n)}m
=
\begin{cases}
x_\alpha^{(n)}m
&\text{si $p|n$}
\\
0
&\text{sinon}.
\end{cases}
\end{align*}
De m\^eme avec $x_{-\alpha}^{(n)}$,
et donc
\begin{equation}
x_{\pm\alpha}^{(n)} f_m=
\begin{cases}
 f_{x_{\pm\alpha}^{(n)}m}&\text{si $p|n$}
\\
0
&\text{sinon}.
\end{cases}
\end{equation}

D'autre part,
pour
$m\in(M^\phi)^{[1]}$,
on a
\begin{align*}
x_{\pm\alpha}^{(n)}\bullet m
&=
\Dist(F)(x_{\pm\alpha}^{(n)})\bullet m
\\
&\hspace{1,5cm}
\quad\text{consid\'erant $m$ comme un \'el\'ement de $M^\phi$ dans le membre de droite}
\\
&=
\begin{cases}
x_{\pm\alpha}^{(\frac{n}{p})}\bullet
m
&\text{si $p|n$,
}
\\
0
&\text{sinon}
\end{cases}
\\
&=
\begin{cases}
x_{\pm\alpha}^{(n)}m
&\text{si $p|n$,
consid\'erant cette fois $m$
comme un \'el\'ement de
$M$,}
\\
0
&\text{sinon,}
\end{cases}
\end{align*}
et donc
(3) devient
$x_{\pm\alpha}^{(n)} f_m=f_{x_{\pm\alpha}^{(n)}\bullet m}$, et l'\'egalit\'e (2) est donc vraie, comme voulu.
\\

\setcounter{equation}{0}
\noindent
(2.4)
Le foncteur de  torsion par Frobenius des    $G$-modules \'etant un adjoint \`a droite du foncteur
$(?^{G_1})^{[-1]}$, on a donc
\begin{cor}
(i)
Le foncteur de contraction par Frobenius est un adjoint \`a droite du foncteur $\St\otimes\St\otimes(?^{[1]})$
sur la cat\'egorie des $G$-modules.

(ii)
Si $M$ est un  $G$-module injectif,  $M^\phi$ l'est \'egalement.

\end{cor}

\pf
(i)
Pour deux
$G$-modules
$M$ et $V$,
on a
\begin{align*}
G\Mod
&
(\St\otimes\St\otimes
V^{[1]}, M)
\simeq
G\Mod(V^{[1]},\St\otimes\St\otimes
M)
\\
&\simeq
G\Mod(V^{[1]},
(\St\otimes\St\otimes
M)^{G_1})
\\
&\simeq
G\Mod(V^{[1]},
(((\St\otimes\St\otimes
M)^{G_1})^\phi)^{[1]})
\\
&\simeq
G\Mod(V,
((\St\otimes\St\otimes
M)^{G_1})^\phi))
\\
&\hspace{3cm}
\quad\text{car
le foncteur   torsion par Frobenius est pleinement fid\`ele}
\\
&\simeq
G\Mod(V,
M^\phi).
\end{align*} 

(ii)
C'est imm\'ediat \`a partir de  (i) ; en effet, 
$\St\otimes\St\otimes(?^{[1]})$ est un foncteur  exact.

\setcounter{equation}{0}
\noindent
(2.5)
Soit $r\in\bbN^+$,
$\St_r=L((p^r-1)\rho)$ le $r$-i\`eme module de Steinberg.
It\'erant le foncteur de contraction par Frobenius, on obtient
\begin{thm}
La $r$-i\`eme   contraction par Frobenius
$(?)^{\phi^r}$
est adjointe \`a droite du foncteur
$\St_r\otimes\St_r\otimes(?^{[r]})$.

\end{thm}

\setcounter{equation}{0}
\begin{center}
$3^\circ$
{\bf 
Bonnes filtrations}

\end{center}

Pour tout $\lambda\in\Lambda^+$, rappelons qu'on dispose des \emph{modules standards} ou \emph{modules de Weyl duaux} $\nabla(\lambda)=\{f\in\Sch_\Bbbk(G,\lambda)|f(gb)=b^{-1}f(g),\ \forall g\in G, \ \forall b\in B\}$. On peut les d\'efinir sur    $\bbZ$, et les $\nabla_\bbZ(\lambda)\otimes_\bbZ\bbQ$ fournissent (avec les notations \'evidentes) les modules simples pour
$G_\bbZ\otimes_\bbZ\bbQ$.
On a, par exemple, 
$\St=\nabla((p-1)\rho)$. On dispose \'egalement, pour tout  $\lambda\in\Lambda^+$, des   \emph{modules de Weyl}  $\Delta(\lambda)$
de plus haut poids $\lambda\in\Lambda^+$
\cite[II.2.13 Rmk 1)]{J}. 

Nous dirons qu'un  $G$-module de dimension finie $M$ admet une \emph{filtration de    Weyl}  si et seulement si tous les sous-quotients correspondants sont   des modules de Weyl. 
Dualement, on dira que
$M$ admet une \emph{bonne filtration} si et seulement si  tous les sous-quotients correspondants sont de la forme   
$\nabla(\lambda)$.

\noindent
(3.1) On va s'int\'eresser \`a l'existence de bonnes filtrations sur $M^\phi$.
Rappelons qu'il existe 
\cite[\S3 ]{vdK}
un $G$-module $M$  
admettant une bonne filtration tel que  
$(M^{G_1})^{[-1]}$ n'en admette pas.
N\'eanmmoins,  pour tout $G$-module, on a  
$
M^\phi\geq
(M^{G_1})^\phi=
\{\{(M^{G_1})^{[-1]}\}^{[1]}\}^\phi=
(M^{G_1})^{[-1]}$. Rappelons aussi que l'application naturelle  
$\nabla(p\lambda)^\phi \to \nabla(\lambda)$ ($\lambda\in\Lambda^+$) est naturellement scind\'ee par l'application d'\'el\'evation \`a la puissance $p$-i\`eme \cite[Prop. 3.2]{GK11}, faisant donc de $\nabla(\lambda)$ un facteur direct de $\nabla(p\lambda)^\phi$.

Pour \'enoncer plus simplement les r\'esultats qui vont suivre, nous adopterons par commodit\'e  la convention 
que le module nul admet \`a la fois une filtration de   Weyl et une bonne filtration 
(l'ensemble de ces filtrations \'etant toutefois r\'eduit \`a l'ensemble vide). Le r\'esultat suivant am\'eliore inconditionellement \cite[Cor. 5.8.1]{GK}

\begin{thm}
Soit  $M$un $G$-module de dimension finie.
Si $M$ admet une bonne filtration, il en est de m\^eme de
$M^\phi$.
\end{thm}

\pf
Il suffit de montrer que  
$\Ext^1(\Delta(\lambda),M^\phi)=0$ 
pour tout $\lambda\in\Lambda^+$
\cite[II.4.16]{J}.
Mais
\begin{align*}
\Ext^1
(\Delta(\lambda),M^\phi)
&\simeq
\Ext^1(\St\otimes\St\otimes\Delta(\lambda)^{[1]},M)
\\
&\hspace{3cm}
\text{car le foncteur de contraction par  Frobenius est
exact,
}
\\
&\simeq
\Ext^1(\St\otimes\Delta((p-1)\rho+p\lambda),M)
\\
&\hspace{0,5cm}
\text{par la version duale d'un th\'eor\`eme de  
Andersen-Haboush \cite[II.3.19]{J}}
\\
&=0
\quad\text{car $\St\otimes\Delta((p-1)\rho+p\lambda)$
admet une filtration de Weyl
\cite[II.4.21]{J}}.
\end{align*}

\setcounter{equation}{0}
\noindent
(3.2)
{\bf Remarque.}
Utilisant la 
$G_1$-projectivit\'e du module de  Steinberg  
\cite[II.10.2]{J}, Donkin r\'eduit le probl\`eme \`a montrer que  
$G_1\Mod(\St,\nabla(\lambda))^{[-1]}$ reste un   $G$-module standard s'il ne s'annule pas.
S'il est non nul,
le $G_1$-socle de $\nabla(\lambda)$ doit \^etre dans le m\^eme bloc que celui de la repr\'esentation de  Steinberg, 
et donc
$\lambda=(p-1)\rho+p\lambda^1$ pour un certain 
$\lambda^1\in\Lambda^+$.
On a alors
$\nabla(\lambda)\simeq\St\otimes\nabla(\lambda^1)^{[1]}$
de nouveau gr\^ace au th\'eor\`eme  d'Andersen-Haboush {\cite[II.3.19]{J}}, et
$G_1\Mod(\St,\nabla(\lambda))^{[-1]}\simeq
\nabla(\lambda^1)$
par $G_1$-semi-simplicit\'e de $\St$.

\setcounter{equation}{0}
\noindent
(3.3)
{\bf Corollaire.}
{\it
(i)
Pour un  $G$-module de dimension finie
$M$
admettant une bonne filtration, notons 
$(M:\nabla(\lambda))$ la multiplicit\'e de 
$\nabla(\lambda)$,
$\lambda\in\Lambda^+$ dans celle-ci. Alors, pour tous 
$\lambda,\mu\in\Lambda^+$, on a 
\begin{align*}
(\nabla(\lambda)^\phi:\nabla(\mu))
&=
(\St\otimes\nabla((p-1)\rho+p\mu):\nabla(\lambda))
=
(\St\otimes\nabla(\lambda):
\nabla((p-1)\rho+p\mu))
\\
&=
\sum_{w\in W/\rC_{W\bullet}(\mu)}(-1)^{\ell(w)}\dim\nabla(\lambda)_{p(w\bullet\mu)}
\\
&=
\#\{\pi\in\bbB(\lambda)|\text{$\pi$ {\rm{est}}\, $(p-1)\rho$-{\rm{dominant et}} $\pi(1)=p\mu$}\} 
\end{align*} 
avec $\rC_{W\bullet}(\mu)=\{w\in W|w\bullet\mu=\mu\}$
et
$\bbB(\lambda)$ est l'ensemble des  chemins de Lakshmibai-Seshadri de  forme $\lambda$
{\rm{\cite[2, Def]{Li94}}}.
En particulier, si $\mu\in\Lambda^+$ est maximal parmi ceux  des 
$\nu\in\Lambda$ tels que 
$p\nu$ soit un poids de  $\nabla(\lambda)$, on a
\[
(\St\otimes\nabla(\lambda):
\nabla((p-1)\rho+p\mu))
=
\dim\nabla(\lambda)_{p\mu}.
\]

(ii)
Un $G$-module de basculement (``tilting module'' {\rm{\cite[E.1]{J}}})  admet une contraction par Frobenius qui est  $G$-module de basculement.

}

\pf
(i)
On a 
\begin{align*}
(\nabla(\lambda)^\phi
&
:\nabla(\mu))
=
\dim
G\Mod(\Delta(\mu),\nabla(\lambda)^\phi)
\quad\text{gr\^ace \`a \cite[II.4.16]{J}}
\\
&=
\dim
G\Mod(\St\otimes\St\otimes\Delta(\mu)^{[1]},\nabla(\lambda))
\\
&=
\dim
G\Mod(\St\otimes\Delta((p-1)\rho+p\mu),\nabla(\lambda))
=
(\St\otimes\nabla((p-1)\rho+p\mu):\nabla(\lambda))
\\
&=
\dim
G\Mod(\Delta((p-1)\rho+p\mu),\St\otimes\nabla(\lambda))
=
(\St\otimes\nabla(\lambda):\nabla((p-1)\rho+p\mu)).
\end{align*}
Pour un  $B$-module $M$, introduisons
 $\chi(M)=\sum_{i\in\bbN}(-1)^i\ch\rR^i\ind_B^G(M)$ son caract\`ere d'Euler  \cite[II.5.7]{J}, ch d\'esignant le caract\`ere formel \cite[I.2.11(6)]{J}.
Alors
\begin{align*}
\ch(\St\otimes\nabla(\lambda))
&=
\ch\nabla((p-1)\rho\otimes\nabla(\lambda))
\quad\text{gr\^ace \`a l'identit\'e tensorielle {\rm{\cite[I.3.6]{J}}}}
\\
&=
\chi((p-1)\rho\otimes\nabla(\lambda))
\\
&\hspace{3cm}
\quad\text{gr\^ace au th\'eor\`eme d'annulation de Kempf {\rm{\cite[II.4.6]{J}}}}
\\
&=
\sum_{\nu\in\Lambda}\chi((p-1)\rho+\nu)\dim\nabla(\lambda)_\nu
\quad\text{car
$\chi$ est additive}.
\end{align*}
Pour $w\in W$
on a
$\chi(w^{-1}\bullet((p-1)\rho+\nu))
=
(-1)^{\ell(w)}\chi((p-1)\rho+\nu)
$,
et
$w^{-1}\bullet((p-1)\rho+\nu)=(p-1)\rho+p\mu$
si et seulement si
$\nu=p(w\bullet\mu)$.
Les composants  de 
$\ch(\St\otimes\nabla(\lambda))$
contribuant \`a la multiplicit\'e de
$(\St\otimes\nabla(\lambda):\nabla((p-1)\rho+p\mu))$ sont donc juste les
\[
\sum_{w\in W/\rC_{W\bullet}(\mu)}
\chi(w\bullet((p-1)\rho+p\mu))
\dim\nabla(\lambda)_{p(w\bullet\mu)},
\]
et donc
$(\St\otimes\nabla(\lambda):\nabla((p-1)\rho+p\mu))=
\sum_{w\in W/\rC_{W\bullet}(\mu)}
(-1)^{\ell(w)}
\dim\nabla(\lambda)_{p(w\bullet\mu)}$.

(ii) d\'ecoule du th\'eor\`eme et de  (1.2.v).

\setcounter{equation}{0}
\noindent
(3.4) Montrons maintenant la n\'ecessit\'e des hypoth\`eses de caract\'eristique puis de petitesse du plus haut poids dans la proposition 1.3 sur l'exemple d'un groupe $G$ de type $\rG_2$ dont nous noterons $\alpha$ et $\beta$ les racines simples (avec $\alpha$ la courte)   et  
$\varpi_\alpha$ et $\varpi_\beta$ les poids fondamentaux correspondants.

(i) Supposons que  $p=2$. Les poids de  $\St=L(\rho)$ appartenant \`a $2\Lambda^+$ sont juste  $2\varpi_\alpha$ et $0$ de multiplicit\'es respectivement 2 et 4.
Par (3.1),
$\St^\phi=\nabla(\rho)^\phi$ admet une bonne filtration avec comme sous-quotients 
$\nabla(\varpi_\alpha)$ et $\Bbbk$ apparaissant chacun deux fois.
Cependant $\nabla(\varpi_\alpha)$ n'est pas simple, ayant 
$\Bbbk$ comme module de t\^ete par la formule des sommes de Jantzen 
\cite[II.8.19]{J}.
Il s'ensuit que  $L(\rho)^\phi$ ne peut \^etre  semi-simple dans ce cas. 

(ii) Revenons à la situation  et notations   g\'en\'erales  de (1.3). M\^eme si
 $L(\lambda^0)^\phi$ reste semi-simple, 
$L(\lambda^0)^\phi\otimes
L(\lambda^1)$ n'a aucune raison de l'\^etre.
En effet, consid\'erons
$L(\lambda)\otimes\St^{[1]}$ avec $\lambda=(p-3)\varpi_\alpha+2\varpi_\beta$ en caract\'eristique    $p\geq3$.
Comme $p\varpi_\alpha=\lambda-\alpha_2$ est un poids de $\nabla(\lambda)$,
c'en est \'egalement un  de $L(\lambda)$ ; si $L(\mu)$ d\'esigne un autre facteur de composition de $\nabla(\lambda)$,  alors
$\lambda-\alpha_2\not\leq\mu$.
Ainsi,
$\varpi_\alpha$ est un poids maximal de
$L(\lambda)^\phi$, et par suite
$L(\lambda)^\phi$ a un facteur de composition $L(\varpi_\alpha)$.
Comme $L(\varpi_\alpha)=\nabla(\varpi_\alpha)$ pour $p\ne2$,
$(L(\lambda)\otimes\St^{[1]})^\phi\simeq
L(\lambda)^\phi\otimes\St$ admet un  quotient $\nabla(p\varpi_\alpha+(p-1)\varpi_\beta)$.
Mais
$\nabla(p\varpi_\alpha+(p-1)\varpi_\beta)$ a un facteur de composition
$L(p\varpi_\alpha+(p-1)\varpi_\beta-\alpha)$,
et n'est donc pas  semi-simple.

\setcounter{equation}{0}
\begin{center}
$4^\circ$
{\bf 
$G_r$-modules}

\end{center}

Soit $r\in\bbN$ avec $r\geq2$.

\setcounter{equation}{0}
\noindent
(4.1)
La contraction par  Frobenius   induit un foncteur de la  cat\'egorie des $G_r$-modules dans celle des $G_{r-1}$-modules,
que l'on continue de noter    $?^\phi$.

\begin{thm}
Il existe un isomorphisme de foncteurs  
$?^\phi\to
\{(\St\otimes\St\otimes?)^{G_1}\}^{[-1]}$
de la  cat\'egorie des $G_r$-modules dans celle des $G_{r-1}$-modules.
\end{thm}

\pf
Soit
$M$ un $G_r$-module.
On montre que
$M^\phi\to
G_1\Mod(\St\otimes\St, M)^{[-1]}$
via
$m\mapsto f_m$
est
$G_{r-1}$-lin\'eaire.
Comme $G_{r-1}$ n'est pas, en g\'en\'eral,  engendr\'e par les
$x_{\pm\alpha}^{(n)}$, $\alpha\in R^\rs$, $n\in]0,p^{r-1}[$,
l'argument pour les  $G$-modules ne se d\'ecalque pas tel quel.
N\'eanmmoins, on peut supposer   $M$ de dimension finie.
Alors, comme les deux foncteurs de l'\'enonc\'e sont exacts,
on a seulement \`a v\'erifier la  $G_{r-1}$-lin\'earit\'e de l'application ci-dessus sur les   $G_r$- modules simples.
Mais ces derniers sont \'equip\'es d'une structure de $G$-modules, et l'on sait d\'ej\`a que l'application  
est  $G$-lin\'eaire par  (2.3) et donc,  a fortiori,
$G_{r-1}$-lin\'eaire.

\setcounter{equation}{0}
\noindent
(4.2)
De m\^eme pour les  $G_rT$-modules et pour les  $G_rB$-modules.

\begin{cor}
Il existe un isomorphisme de foncteurs
$?^\phi\to
\{(\St\otimes\St\otimes?)^{G_1}\}^{[-1]}$
de la cat\'egorie des $G_rT$-(resp. $G_rB$-)
modules dans celle des $G_{r-1}T$-
(resp. $G_rB$-)
modules.

\end{cor}

\pf
Soit 
$M$ un $G_rT$-module.
L'application $m\mapsto f_m$ est $T$-lin\'eaire ;
si $m\in M_{p\lambda}$,
$t\bullet m=\lambda(t)m$.
Regardant $f_m$ comme un \'el\'ement de
$G_1\Mod(\St\otimes\St,M)$,
$tf_m=tf_m(t^{-1}?)$.
Comme $tf_m(t^{-1}(v_+\otimes v_-))=
tm=(p\lambda)(t)m$,
$tf_m=(p\lambda)(t)f_m$, et donc
$f_m$
a comme poids
$\lambda$
comme \'el\'ement de
$G_1\Mod(\St\otimes\St,M)^{[-1]}$. 
Comme l'application ci-dessus est  $G_{r-1}$-lin\'eaire gr\^ace au th\'eor\`eme,
l'assertion pour les $G_rT$-modules s'ensuit.
En ce qui concerne les $G_rB$-modules, 
$\Dist(U)$ est engendr\'e par
$x_{-\alpha}^{(n)}$, $\alpha\in R^\rs$, $n\in\bbN$,
et donc l'application est  $\Dist(U)$-lin\'eaire par  (2.3). L'assertion s'ensuit.

\setcounter{equation}{0}
\noindent
(4.3)
Il d\'ecoule de tout cela que la contraction par Frobenius sur les   categories des $G_r$-modules,
$G_rT$-modules,
$G_rB$-modules
sont adjoints \`a droite de la torsion par Frobenius tensoris\'e deux fois avec le module de Steinberg.  
Plus g\'en\'eralement,

\begin{cor}
Soit $s\in\bbN$ avec $0<s<r$.

(i)
Le foncteur de  $s$-i\`eme contraction par  Frobenius 
$?^{\phi^s}$de la cat\'egorie des
$G_r$-(resp. $G_rT$-, $G_rB$-)
modules dans celle des $G_{r-s}$-
(resp. $G_{r-s}T$-, $G_{r-s}B$-)
modules
est adjoint \`a droite du foncteur
$\St_s\otimes\St_s\otimes(?^{[s]})$ de la cat\'egorie des
$G_{r-s}$-(resp. $G_{r-s}T$-, $G_{r-s}B$-)
modules dans celle des $G_{r}$-
(resp. $G_{r}T$-, $G_rB$-)
modules.

(ii)
Les foncteurs 
$?^{\phi^s}$ sur ces cat\'egories 
pr\'eservent l'injectivit\'e et la projectivit\'e des modules de dimension finie.

\end{cor}

\setcounter{equation}{0}
\noindent
(4.4)
Pour un 
$B_r$-
(resp. $B_rT$-, $B$-)
module $M$
posons
$\nabla_r(M)=\ind_{B_r}^{G_r}(M)$
(resp.
$\hat\nabla_r(M)=\ind_{B_rT}^{G_rT}(M)$,
$\tilde\nabla_r(M)=\ind_{B}^{G_rB}(M)$).
On dira qu'un  $G_r$-
(resp.
$G_rT$-,
$G_rB$-)
module $V$ de dimension finie admet une
$\nabla_r$-
(resp.
$\hat\nabla_r$-,
$\tilde\nabla_r$-)
\emph{filtration} si et seulement si  $V$ admet une  $G_r$-
(resp.
$G_rT$-,
$G_rB$-)
filtration dont tous les 
sous-quotients correspondants sont de la forme
$\nabla_r(\lambda)$
(resp.
$\hat\nabla_r(\lambda)$,
$\tilde\nabla_r(\lambda)$),
$\lambda\in\Lambda$.  
Pour un  $B_r^+T$-module $M'$, posons aussi
$\hat\Delta_r(M')=\coind_{B_r^+T}^{G_r^+T}(M')$
\cite[II.9.1.5]{J}.

Rappelons  (\cite[II.11.2]{J}) qu'un   
$G_rT$-module $M$ de dimension finie admet une  $\hat\nabla_r$-filtration si et seulement si
$M$ est $B^+_rT$-injectif.
Comme pour les $G$-modules, on a 

\begin{lem}
Un $G_rT$-module $M$ de dimension finie admet une  $\hat\nabla_r$-filtration si et seulement si
$\Ext^1_{G_rT}(\hat\Delta_r(\lambda),M)=0$ pour tout 
$\lambda\in\Lambda$.

\end{lem}

\pf
Supposons tout d'abord que $\Ext_{G_1T}^1(\hat\Delta(\lambda),M)=0$
pour tout $\lambda\in\Lambda$.
Si $Q$ est une 
$B_1^+T$-enveloppe injective du $B_1^+T$-socle de$M$,
on veut voir que le monomorphisme 
$M\hookrightarrow Q$ est scind\'e,
ce qui d\'ecoulera de 
$\Ext_{B_1^+T}^1(Q/M,M)=0$.
Pour avoir cela, il suffit de v\'erifier que
$\Ext_{B_1^+T}^1(\lambda,M)=0$
pour tout $\lambda\in\Lambda$.
Mais
\begin{align*}
\Ext_{B_1^+T}^1(\lambda,M)
&\simeq
\Ext_{B_1^+T}^1(M^*,-\lambda)\simeq
\Ext_{G_1T}^1(M^*,\ind_{B_1^+T}^{G_1T}(-\lambda))
\\
&\simeq
\Ext_{G_1T}^1(M^*,\hat\Delta(-\lambda+2(p-1)\rho))
\quad\text{gr\^ace \`a \cite[II.9.2]{J}}
\\
&\simeq
\Ext_{G_1T}^1(\hat\Delta(-\lambda+2(p-1)\rho)^*,M)
\\
&\simeq
\Ext_{G_1T}^1(\hat\Delta(\lambda),M)
\quad\text{de nouveau gr\^ace \`a \cite[II.9.2]{J}}
\\
&=0.
\end{align*}
La r\'eciproque vaut \'egalement.

\setcounter{equation}{0}
\noindent
(4.5)
Soit $s\in\bbN$ avec $0<s<r$.
\begin{thm}
Si $M$ est un  $G_rT$-module de dimension finie admettant une   $\hat\nabla_r$-filtration, 
$M^{\phi^s}$ admet une  $\hat\nabla_{r-s}$-filtration.

\end{thm}

\pf
Par (4.4) il suffit simplement de voir que
$\Ext^1_{G_{r-s}T}(\hat\Delta_{r-s}(\lambda),M^{\phi^s})=0$
pour tout $\lambda\in\Lambda$.
Par (4.3) 
$\Ext^1_{G_{r-s}T}(\hat\Delta_{r-s}(\lambda),M^{\phi^s})\simeq
\Ext^1_{G_{r}T}(\St_s\otimes\St_s\otimes\hat\Delta_{r-s}(\lambda)^{[s]},M)$.
Maintenant
$\hat\Delta_{r-s}(\lambda)$ s'\'etend en
$\tilde\Delta_{r-s}(\lambda)$
\cite[II.9.1.4]{J},
et, comme $G$-modules,
\begin{align*}
\St_s
\otimes\tilde\Delta_{r-s}(\lambda)
&\simeq
\St_s\otimes\ind_{B^+}^{G_{r-s}B^+}(\lambda-2(p^{r-s}-1)\rho)^{[s]}
\quad\text{gr\^ace \`a
\cite[II.9.2]{J}}
\\
&\simeq
\St_s\otimes\ind_{G_sB^+/G_s}^{G_{r}B^+/G_s}(p^s(\lambda-2(p^{r-s}-1)\rho))
\\
&\simeq
\St_s\otimes\ind_{G_sB^+}^{G_{r}B^+}(p^s(\lambda-2(p^{r-s}-1)\rho))
\quad\text{gr\^ace \`a
\cite[II.3.19]{J}}
\\
&\simeq
\ind_{G_sB^+}^{G_{r}B^+}(\St_s\otimes
p^s(\lambda-2(p^{r-s}-1)\rho))
\\
&\hspace{6,5cm}
\quad\text{gr\^ace \`a l'identit\'e tensorielle
\cite[I.3.6]{J}}
\\
&\simeq
\ind_{G_sB^+}^{G_{r}B^+}(\ind_{B^+}^{G_sB^+}((p^s-1)\rho-2(p^s-1)\rho)\otimes
p^s(\lambda-2(p^{r-s}-1)\rho))
\\
&\hspace{8cm}
\text{ de nouveau gr\^ace \`a
\cite[II.9.2]{J}}
\\
&\simeq
\ind_{G_sB^+}^{G_{r}B^+}(\ind_{B^+}^{G_sB^+}(-(p^s-1)\rho+
p^s(\lambda-2(p^{r-s}-1)\rho)))
\\
&\simeq
\ind_{B^+}^{G_{r}B^+}((p^s-1)\rho+
p^s\lambda-2(p^r-1)\rho)
\\
&\simeq
\hat\Delta_r((p^s-1)\rho+p^s\lambda)
\quad\text{gr\^ace \`a
\cite[II.9.2]{J}}.
\end{align*}
Ainsi
\begin{align*}
\St_s\otimes\St_s\otimes\hat\Delta_{r-s}(\lambda)^{[s]}
&\simeq
\St_s\otimes
\hat\Delta_{r}((p^s-1)\rho+p^s\lambda)
\\
&\simeq
\hat\Delta_{r}(\St_s\otimes((p^s-1)\rho+p^s\lambda))
\\
&\hspace{5cm}
\text{gr\^ace \`a l'identit\'e tensorielle \cite[I.3.6]{J}}
\end{align*}
admet une $\hat\Delta_r$-filtration, et donc
$\Ext^1_{G_{r}T}(\St_s\otimes\St_s\otimes\hat\Delta_{r-s}(\lambda)^{[s]},M)=0$ par
(4.4), comme voulu.

\setcounter{equation}{0}
\noindent
(4.6)
Soit $s\in\bbN$ avec $0<s<r$.
\begin{cor}
Si $M$  est un  $G_r$-
(resp. $G_rB$-)
module  de dimension finie admettant une $\nabla_r$-
(resp.
$\tilde\nabla_r$-)
filtration, 
$M^{\phi^s}$ admet une $\nabla_{r-s}$-
(resp.
$\tilde\nabla_{r-s}$-)
filtration.

\end{cor}

\pf
Pour les $G_r$-modules, il suffit de voir que chaque
$\nabla_r(\lambda)^{\phi^s}$,
$\lambda\in\Lambda$,
admet une  $\nabla_{r-s}$-filtration.
Mais
$\nabla_r(\lambda)$ s'\'etend \`a
$\hat\nabla_r(\lambda)$, et l'assertion d\'ecoule donc de  (4.5).

Pour les $G_rB$-modules
il suffit de nouveau de voir que chaque 
$\tilde\nabla_r(\lambda)^{\phi^s}$,
$\lambda\in\Lambda$,
admet une $\tilde\nabla_{r-s}$-filtration.
Comme $\tilde\nabla_r(\lambda)=
\hat\nabla_r(\lambda)$ en tant que $G_rT$-modules,
$\tilde\nabla_r(\lambda)^{\phi^s}$
admet une $\hat\nabla_{r-s}$-filtration par (4.5), et donc
$\tilde\nabla_r(\lambda)=
\hat\nabla_r(\lambda)$ est $B_{r-s}^+T$-injectif par
(4.4).
Le module
$\tilde\nabla_r(\lambda)^{\phi^s}$
admet donc une $\tilde\nabla_{r-s}$-filtration par
\cite[Rmk. II.11.2.2]{J}.

\end{document}